\documentclass[12pt,reqno]{amsart}
\usepackage{amssymb,latexsym,amsmath}
\usepackage{setspace}
\newtheorem{theorem}{Theorem}

\theoremstyle{definition}
\newtheorem{definition}{Definition}

\begin{document}
\begin{center}\textbf{BesselK Derivatives with respect to Order at one half}\end{center}
\begin{center}\textbf{Charles Ryavec}\end{center}
The order derivatives of the modified Bessel function of the second kind at $s = .5$ are obtained as finite expressions of integrals that generalize the exponential integral appearing in the first derivative (Theorem 1). These derivatives arise in the investigation of a BesselK relationship with the Riemann zeta function. Any use of the term, derivative, with respect to a Bessel function, here means a derivative with respect to its order.
\section {Introduction}
Given multiplicative arithmetic functions,
\begin{eqnarray*}
a(n) &=& \sigma_1(n) \qquad (n, 2) = 1\\
&=& 2^m \qquad \qquad n = 2^m \\ 
b(n) &=& \sigma_1(n) \qquad (n, 2) = 1\\
&=& 0 \qquad \qquad n = 2^m, \ m > 0, 
\end{eqnarray*}
for each positive integer, $j$, define a function of the complex variable, $s$, by
$$
c[s, j] = \sum_{d|j} a(d) \ b(\frac{j}{d}) \ \Big(\frac{j}{d^2}\Big)^{\frac{s}{2}}.
$$
Define an entire function of $s$ by,
$$
h(s) = \sum_{j=1}^{\infty} c[s, j] K[s, 2 \pi \sqrt{j}],
$$
where $K[s,x]$ is the modified Bessel function of the second kind. It was shown in ([1]) that,  
$$
h(s) = \frac{s(s+1)}{32\pi^2 \sqrt{2}}(2^{\frac{s}{2}} - 2^{-\frac{s}{2}})(2^{\frac{s-1}{2}}-2^{-\frac{(s-1)}{2}})\zeta^{*}(s)\zeta^{*}(s+1).
$$
Among the zeros of $h(s)$ are the zeros of 
$$
\zeta^{*}(s) = \pi^{-\frac{s}{2}}\Gamma(\frac{s}{2})\zeta(s).
$$.  
\section {Motivation}
The dual representation of $h(s)$ provided the motivation to study $K[s, x]$. The special interest attached to the $\sigma = .5$ line raised the consideration of the coefficients in the expansion,
$$
h(s) = \sum_{n=0}^{\infty}\alpha_n \frac{(s - .5)^n}{n!}
$$
at $s = .5$, one of the special points of, $K[s, x]$, where the value, 
$$
K[.5, x] = \sqrt{\frac{\pi}{2x}} e^{-x},
$$
is simple. The fact ([2], 10.38.7) that the first derivative,
$$
\frac{d}{ds}K[s, x]\Big|_{s = .5} = K[.5, x]\int_0^{\infty}\frac{e^{-u}}{u+2x}du ,
$$
is also simple, with a correspondingly rapid rate of decrease, raised the question of the structure of all the derivatives,
$$
\frac{d^n}{ds^n} K[s, x]\Big|_{s = .5}.
$$
This question combined with another fact, that for each $s$, the function $c[s, j]$ is multiplicative in $j$, so there is a factorization,
$$
c[s, j]  = \prod_{p^e||j}c[s, p^e] . 
$$
Expanding,  
$$
c[s, p^e] = \sum_{n=0}^{\infty}c_n(p^e)\frac{ (s - .5)^n}{n!}
$$
and
$$
K[s, 2\pi\sqrt{j}] = \sum_{n=0}^{\infty}C_n(j) \frac{(s - .5)^n}{n!},
$$
then
$$
\prod_{p^e||j}\sum_{n=0}^{\infty}c_n(p^e)\frac{(s- .5)^n}{n!}  \sum_{n=0}^{\infty}C_n(j)\frac{(s- .5)^n}{n!} = \sum_{n=0}^{\infty}\alpha_{n}( j)\frac{(s- .5)^n}{n!},
$$
gives coefficients, $\alpha_{n}[j ]$, as polynomial expressions in primes and the derivatives of $K[s, 2\pi\sqrt{j}]$ at $s=.5$. For example, if $j = p$, is a prime, then
$$
c_n(p) = (p + 1)\frac{d^n}{ds^n}(p^{\frac{s}{2}} + p^{-\frac{s}{2}})\Big|_{s=.5} \qquad p > 2
$$
$$
c_n(2^e) = \frac{d^n}{ds^n}2^{e(1 - \frac{s}{2})}\Big|_{s=.5}
$$
$$
\alpha_n(p) = \sum_{m=0}^n {n\choose m}c_m(p)C_{n-m}(p).
$$
The explicit form of the derivatives at $s= .5$ given in Theorem 2 provides the factor, $e^{-2\pi\sqrt{j}}$, in the $C_n(j)$, and they all diminish sufficiently rapidly that the sum on $j$ and $n$ can be reversed. Thus, 
$$
h(s) = \sum_{n=0}^{\infty} \alpha_n \frac{(s-.5)^n}{n!},
$$
$$
\alpha_n = \sum_{j=1}^{\infty}\alpha_n(j).
$$
\section {BesselK Representation}
There are a half dozen basic expressions for $K[s, x]$. Consequently, different expressions for the derivatives exist. A search of the literature did not locate results for derivatives of order greater than the first, even for the special value, $s = .5$, but results for higher derivatives might exist. The fact that relations between primes and order derivatives of K control a result as elusive as the Riemann hypothesis was the motivation to find manageable expressions for the derivatives. The very simple appearance of the first derivative in  ([2], 10.38.7; 6.2.2) suggested 
consideration of the representation,
$$
K[s, x] = \sqrt{\frac{1}{\pi}}\Big(\frac{2}{x}\Big)^s\Gamma[s + \frac{1}{2}]\int_0^{\infty} \frac{\cos[x u]}{(u^2+1)^{s+\frac{1}{2}}}du, \qquad \sigma > -\frac{1}{2}, x > 0.
$$
Define,
$$
T[s, x] = \int_0^{\infty} \frac{\cos[x u]}{(u^2+1)^{s+\frac{1}{2}}}du \qquad \sigma > -\frac{1}{2}, x > 0.
$$
We first obtain the derivatives of  $T[s,x]$. To get the $K$ derivatives, just include the factor,
$$
\sqrt{\frac{1}{\pi}}\Big(\frac{2}{x}\Big)^s\Gamma[s + \frac{1}{2}],
$$
and the additional derivatives. It is instructive to see the derivation of the derivatives,
$$
\frac{d^n}{ds^n}T[s, x],
$$
first for $n = 0, 1$ and then the general case.
\section {Theorem 1}
\begin{theorem}Let $T[s, x]$ be defined for $\sigma > .5, x > 0$, as,
$$
T[s, x] = \int_0^{\infty}\frac{cos[ x u]}{(u^2 + 1)^{s + \frac{1}{2}}}du.
$$
Then,
\begin{eqnarray*}
T[.5, x] &=& \frac{\pi}{2} e^{-x} \\
\frac{d}{ds}T[s, x]\Big|_{s= .5} &=& \Big(\int_0^{\infty}e^{-u}\frac{du}{u+2x} +Log[x] - \dot\Gamma[1] - Log[2]\Big) \ T[.5, x]. 
\end{eqnarray*}
\end{theorem}
\begin{proof} The derivation of T[.5, x] is simply to get it on record. By the residue theorem,
\begin{eqnarray*}
T[.5, x] &=& \frac{1}{2} Re\int_{-\infty}^{\infty}\frac{e^{i x u}}{(u^2 + 1)} du \\
&=& \frac{1}{2} Re \frac{\pi}{2\pi i}\int_{-\infty}^{\infty}e^{i x u}\Big(\frac{1}{u - i} - \frac{1}{u + i}\Big) du \\
&=& \frac{\pi}{2} e^{-x}
\end{eqnarray*}
Next, the integrations will be in the upper half cut plane, cut on $[i, i \infty)$. The cut,  $[-i, -i \infty)$, won't enter the calculations.
\begin{eqnarray}
\frac{d}{ds}T[s, x]\Big|_{s= .5} &=& \frac{1}{2} Re\Big(-\int_{-\infty}^{\infty}e^{i x u} \frac{Log[u^2 + 1]}{(u^2 + 1)} du \Big)\\
&=& \frac{1}{2} Re\Big( -\int_{-\infty}^{\infty}e^{i x u}Log[u^2 + 1](\frac{1}{2 i})\Big(\frac{1}{u - i} - \frac{1}{u + i}\Big) du \Big)\\
&=& \frac{1}{2} Re\Big( \frac{1}{2 i}\int_{-\infty}^{\infty}e^{i x u}\frac{Log[u^2 + 1]}{u + i} du \Big) \\
&-& \frac{1}{2} Re\Big( \frac{1}{2 i}\int_{-\infty}^{\infty}e^{i x u}\frac{Log[u^2 + 1]}{u - i} du \Big).
\end{eqnarray}
First,
\begin{eqnarray*}
(3) &=&  \frac{1}{2} Re\Big( \frac{1}{2 i}\int_{-\infty}^{\infty}e^{i x u}\frac{Log[u^2 + 1]}{u + i} du \Big) \qquad u = i y, \ du = i dy\\
&=&  -\frac{1}{2} Re\Big( \frac{1}{2 i}\int_{\infty}^{1}e^{-x y}(\frac{i\pi}{2} + Log[y+1]+ \frac{i\pi}{2} + Log[y-1])\frac{dy}{y+1} \Big) \\
&-&  \frac{1}{2} Re\Big( \frac{1}{2 i}\int_1^{\infty}e^{-x y}(\frac{i\pi}{2} + Log[y+1]+ \frac{-3 i\pi}{2} + Log[y-1])\frac{dy}{y+1} \Big) \\
&=&  \frac{1}{2} Re\Big( \frac{1}{2 i}\int_1^{\infty}e^{-x y}(\frac{i\pi}{2} + Log[y+1]+ \frac{i\pi}{2} + Log[y-1])\frac{dy}{y+1} \Big) \\
&-&  \frac{1}{2} Re\Big( \frac{1}{2 i}\int_1^{\infty}e^{-x y}(\frac{i\pi}{2} + Log[y+1]+ \frac{-3 i\pi}{2} + Log[y-1])\frac{dy}{y+1} \Big) \\
&=&  \frac{1}{2} Re\Big( \frac{1}{2 i}\int_1^{\infty}e^{-x y}(2 \pi i)\frac{dy}{y+1} \Big) \\
&=&  \frac{\pi}{2} \int_1^{\infty}e^{-x y}\frac{dy}{y+1} \\
&=&  \frac{\pi}{2}e^{-x} \int_0^{\infty}e^{-u}\frac{du}{u+2x}  \\
&=&   \int_0^{\infty}e^{-u}\frac{du}{u+2x} \ T[.5, x].
\end{eqnarray*}
Next,
\begin{eqnarray}
(4) &=& -\frac{1}{2} Re\Big( \frac{1}{2 i}\int_{-\infty}^{\infty}e^{i x u}\frac{Log[u^2 + 1]}{u - i} du \Big) \\
 &=&\frac{1}{2} Re\Big( \frac{1}{2 i}\int_{\infty}^{1+\epsilon}e^{-x y}(\frac{i\pi}{2} + Log[y+1]+ \frac{i\pi}{2} + Log[y-1])\frac{dy}{y-1} \Big) \\
  &+&\frac{1}{2} Re\Big( \frac{1}{2 i}\int_{1+\epsilon}^{\infty}e^{-x y}(\frac{i\pi}{2} + Log[y+1]+ \frac{-3 i\pi}{2} + Log[y-1])\frac{dy}{y-1} \Big) \\
  &+& \frac{1}{2} Re\Big( \frac{1}{2 i}\int_{\frac{\pi}{2}}^{\frac{-3\pi}{2}}e^{-x (1 - i\epsilon e^{i\theta})}( Log[2i + \epsilon e^{i\theta}]+ Log[\epsilon e^{i\theta}])i d\theta \Big) \\
  &=&-\frac{\pi}{2} \Big(\int_{1+\epsilon}^{\infty}e^{-x y}\frac{dy}{y-1} \Big) \\
  &+& Re\Big( -\frac{1}{4}\int_{-\pi}^{\pi}e^{-x (1 - i\epsilon e^{i\theta})}( Log[2]+ Log[\epsilon] + i\varphi) d\varphi \Big) \\
  &=&-\frac{\pi}{2} \Big(Log[y-1]e^{- x y}\Big|_{1 + \epsilon}^{\infty} + \int_{1+\epsilon}^{\infty}xLog[y -1]e^{-x y}dy \Big) \\
  &+& Re\Big( -\frac{1}{4}\int_{-\pi}^{\pi}e^{-x }( Log[2]+ Log[\epsilon] + i\varphi) d\varphi \Big) \\
  &=&  -\frac{\pi}{2}\Big(-Log[\epsilon] + \int_0^{\infty}Log[\frac{v}{x}]e^{-v}dv\Big)e^{-x} -\frac{\pi}{2}\Big(Log[2] + Log[\epsilon]\Big)e^{-x}
\end{eqnarray}
where we have let $\epsilon = 0$ in the cases when the limit is obvious. Expression(13) simplifies:
$$
(13) = \Big(Log[x] - \dot\Gamma[1] - Log[2]\Big)T[.5, x].
$$
The final result is, then, (3) + (13), which is 
$$
\frac{d}{ds}T[s, x]\Big|_{s= .5} = \Big(\int_0^{\infty}e^{-u}\frac{du}{u+2x} +Log[x] - \dot\Gamma[1] - Log[2]\Big) \ T[.5, x],  
$$
which is the result stated in the theorem. 
\end{proof}
Note that
\begin{eqnarray*}
\frac{d}{ds}K[s, x]\Big|_{s = .5} &=& \frac{1}{\sqrt{\pi}}\frac{d}{ds}\Big((\frac{2}{x})^s\Gamma[s+ .5]T[s, x]\Big)\Big|_{s=.5} \\
&=& \Big(Log[2] - Log[x] + \dot \Gamma[1] +\frac{\dot T}{T}\Big) K[.5, x] \\
&=& \int_0^{\infty}e^{-u}\frac{du}{u+2x} \ K[.5, x],
\end{eqnarray*}
which is the remarkable logarithmic derivative of $K$. 
\section {Definitions}
Several polynomials will appear in the statement, as well as in the proof, of the next theorem. These polynomials are defined now. Let
\begin{eqnarray*}
U[y] &=& -i\pi + Log[y + 1] + Log[y-1] \\
V[y] &=& i\pi + Log[y + 1] + Log[y-1].
\end{eqnarray*}
Then
$$
U^n - V^n = (U - V)\sum_{k=0}^{n-1}U^{n-1-k}V^k.
$$
\begin{definition}Define polynomials, $p_n$, as follows. The polynomial in variables, $U$ and $V$, 
$$
\sum_{k=0}^{n-1}U^{n-1-k}V^k,
$$
may be put in the variables,
\begin{eqnarray*}  
U + V &=& 2(Log[y + 1] + Log[y-1]) \\
UV &=& \pi^2 +(Log[y + 1] + Log[y-1])^2, 
\end{eqnarray*}
and, so, for $n\ge 1$, it is a polynomial,
$$
p_n(Log[y+1] + Log[y -1] ) 
$$
This is how the polynomial, $p_n$, is defined. \end{definition} For example,
$$
p_2(u) = 2 u.
$$
\begin{definition}
Define polynomials, $f[n, k, y]$, in the variable, $Log[y+1]$, as
$$
f[n, k, y] =  ( i\pi + Log[y+1])^{n-k}- (-i\pi + Log[y+1])^{n-k}
$$
\end{definition}
Then,
$$
f[n, k, 1] =  ( i\pi + Log[2])^{n-k}- (-i\pi + Log[2])^{n-k},
$$
and 
$$
\frac{d}{dy}f[n, k, y] = (n - k)\frac{f[n-1, k, y]}{y+1}. 
$$
\begin{theorem} 
Let $T[s, x], p_n, f[n, k, y]$,  be defined as above. Then
\begin{eqnarray*}
\frac{d^n}{ds^n}T[s, x]\Big|_{s =.5} &=& (-1)^n \int_0^{\infty}cos[ x u]\frac{Log^n[u^2+1]}{u^2 + 1}du \\
&=& (-1)^n T[.5,x]\Big(A_1[n, x] + A_2[n, x] +A_3[n, x] +A_4[n, x] \Big),
\end{eqnarray*}
where,
\begin{eqnarray*}
A_1[n, x] &=& Re\Big(\frac{1}{2\pi}\int_{-\pi}^{\pi}(Log[2] + i\varphi)^n d\varphi\Big) \\
A_2[n, x] &=& -\int_{0}^{\infty} e^{- x u}p_n(Log[u+2]  + Log[u]) \frac{dv}{u + 2} \\
A_3[n, x] &=& x\sum_{k=0}^{n-1}{n\choose k}\int_0^{\infty} e^{-xu}Re\Big(\frac{1}{2\pi i}f[n, k, u+1]\Big)  \frac{Log^{k+1}[u]}{k+1}du \\
A_4[n, x] &=& -\sum_{k=0}^{n-1} {n\choose k}\int_0^{\infty}e^{-xu}Re\Big(\frac{1}{2\pi i}(\frac{d}{du}f[n, k, u+1] )\Big) \frac{Log^{k+1}[u]}{k+1}du  \qquad      
\end{eqnarray*}
\end{theorem}
Postponing the proof momentarily, note that $A_1[n, x]$ is a polynomial in $\pi$ and $Log[2]$, while the expressions, $A_k[n, x],  \ 2\le k\le 4$, consist of linear combinations of integrals of the form,
$$
U[a, b, \epsilon] = \int_0^{\infty} e^{-x u}\frac{Log^{a}[u+2]Log^b[u]}{(u+2)^{\epsilon}} du,
$$
with $a$ and $b$ non negative integers, and $\epsilon = 0, 1$. The case with $a = b = 0$ and $\epsilon = 1$ was encountered in the first derivative (See [2], 6.2.2). Note also the factor, $T[.5, x]$, in all derivatives. The proof will resemble the proof in the $n = 1$ case, with a binomial identity at the end that guarantees that the coefficients of $Log[\epsilon], Log^2[\epsilon], \cdots, Log^n[\epsilon]$, which appear in the calculations, are zero. The limit, $\epsilon \longrightarrow 0$, will give the result. After the proof, the degree of cancellation that occurs in the $n = 2$ case will be considered. 
\begin{proof}From, 
$$
\frac{1}{u^2+1} = \frac{1}{2i}\Big(\frac{1}{u-i} - \frac{1}{u+i}\Big),
$$
we have, 
\begin{eqnarray*}
\frac{d^n}{ds^n}T[s, x]\Big|_{s=.5} &=& (-1)^n \int_0^{\infty}\frac{\cos[x u] Log^n[u^2+1]}{(u^2+1)}du \\
&=& (-1)^n\frac{1}{2} Re\Big(I_1 + I_2  \Big),
\end{eqnarray*} 
$$
I_1 = -\frac{1}{2 i}\int_{-\infty}^{\infty} e^{ix u}\frac{(Log(u + i) + Log[u - i])^n}{(u + i)}du  
$$
$$
I_2 = \frac{1}{2 i}\int_{-\infty}^{\infty} e^{ix u}\frac{(Log(u + i) + Log[u - i])^n}{(u - i)}du.  
$$
Then,
\begin{eqnarray*}
I_1 &=& -\frac{1}{2 i}\int_{-\infty}^{\infty} e^{ix u}\frac{(Log(u + i) + Log[u - i])^n}{(u + i)}du  \\
&=& \frac{1}{2 i}\int_{\infty}^{1} e^{-x y}\Big( \frac{i\pi}{2} + Log[y+1] + \frac{i\pi}{2} + Log [y - 1]\Big)^n\frac{dy}{y + 1} \qquad u = i y  \\
&+& \frac{1}{2 i}\int_1^{\infty} e^{-x y}\Big( \frac{i\pi}{2} + Log[y+1] - \frac{3i\pi}{2} + Log [y - 1]\Big)^n\frac{dy}{y + 1} \\
&=& -\frac{1}{2 i}\int_1^{\infty} e^{-x y}\Big( \frac{i\pi}{2} + Log[y+1] +  \frac{i\pi}{2} + Log[y - 1]\Big)^n\frac{dy}{y + 1} \\
&+& \frac{1}{2 i}\int_1^{\infty} e^{-x y}\Big(  \frac{i\pi}{2} + Log[y+1] - \frac{3i\pi}{2} + Log [y - 1]\Big)^n\frac{dy}{y + 1} \\
&=& \frac{1}{2 i}\int_1^{\infty} e^{-x y}(U^n - V^n)\frac{dy}{y + 1} \\
&=& \frac{1}{2 i}\int_1^{\infty} e^{-x y}(U - V)\sum_{k=0}^{n-1}U^{n-k}V^k\frac{dy}{y + 1}. 
\end{eqnarray*}
Continuing with the definitions from above,
\begin{eqnarray}
I_1 &=& \frac{1}{2 i}\int_1^{\infty} e^{-x y}(U - V)p_n(Log[y+1]+Log[y-1])\frac{dy}{y + 1} \\
 &=& \frac{1}{2 i}\int_1^{\infty} e^{-x y}(-2\pi i)p_n(Log[y+1]  + Log [y - 1])\frac{dy}{y + 1} \\
&=& -\pi \int_{1}^{\infty} e^{-x y}p_n(Log[y+1]  + Log [y - 1])\frac{dy}{y + 1} \\
&=& -\pi e^{-x}\int_{0}^{\infty} e^{- xu}p_n(Log[u+2]  + Log [u])\frac{du}{u + 2} \\
&=& -2T[.5, x] \int_{0}^{\infty} e^{- xu}p_n(Log[u+2]  + Log [u])\frac{du}{u + 2}, 
\end{eqnarray}
an expression again with the factor, $T[.5, x]$. Leave (18) as it is for the moment and do $I_2$ next.
\begin{eqnarray}
I_2 &=& \frac{1}{2 i}\int_{-\infty}^{\infty} e^{ix u}\Big(\frac{Log[u + i] + Log[u - i]\Big)^n}{(u - i)}du  \\
&=& -\frac{1}{2 i}\int_{\infty}^{1 + \epsilon} e^{-x y}\Big(  \frac{i\pi}{2} + Log[y +1] +  \frac{i\pi}{2} + Log[y - 1]\Big)^n\frac{dy}{y - 1}  \\
&-& \frac{1}{2 i}\int_{1 + \epsilon}^{\infty} e^{-x y}\Big(  \frac{i\pi}{2} + Log[y +1] -  \frac{3i\pi}{2} + Log[y - 1]\Big)^n\frac{dy}{y - 1} \\
&-& \frac{1}{2 i}\int_{\frac{\pi}{2}}^{-\frac{3\pi}{2}} e^{-x}\Big(Log[2i] + Log[\epsilon]  + i\theta\Big)^n i d\theta\qquad u = i + \epsilon e^{i\theta}, \ du = i\epsilon e^{i\theta}d\theta \\
&=& -\frac{1}{2i}\int_{1 + \epsilon}^{\infty} e^{-x y}(U^n - V^n)\frac{dy}{y - 1} \\
&+& \frac{1}{2} \int_{\frac{-3\pi}{2}}^{\frac{\pi}{2}} e^{-x}\Big(Log [2] +Log[\epsilon]  + i(\frac{\pi}{2} +\theta)\Big)^n  d\theta,
\end{eqnarray}
where we have let $\epsilon = 0$, in those cases when the limit, $\epsilon\longrightarrow 0$ is obvious. Integrating (23) by parts, there results,
\begin{eqnarray}
-\frac{1}{2i}\int_{1 + \epsilon}^{\infty} e^{-x y}(U^n - V^n )\frac{dy}{y - 1} &=& \frac{1}{2i}\int_{1 + \epsilon}^{\infty} e^{-x y}(V^n - U^n )\frac{dy}{y - 1} \\
&=&
\frac{1}{2i}\int_{1 + \epsilon}^{\infty} e^{-x y}\sum_{k=0}^{n-1} {n\choose k}f[n, k, y]  Log^k[y -1]\frac{dy}{y - 1} \\
&=& \frac{1}{2i}e^{-xy}\sum_{k=0}^{n-1} {n\choose k}f[n, k, y] \frac{Log^{k+1}[y -1]}{k+1}\Big|_{1+\epsilon}^{\infty} \\
&+& \frac{1}{2i}x\int_1^{\infty} e^{-xy}\sum_{k=0}^{n-1}{n\choose k}f[n, k, y]  \frac{Log^{k+1}[y -1]}{k+1}dy \\
&-& \frac{1}{2i}\int_1^{\infty}e^{-xy}\sum_{k=0}^{n-1} {n\choose k}(\frac{d}{dy}f[n, k, y] ) \frac{Log^{k+1}[y -1]}{k+1}dy. 
\end{eqnarray}
The terms (28) and (29) are rewritten as,
$$
(28) = \pi e^{-x}x\sum_{k=0}^{n-1}{n\choose k}\int_0^{\infty} e^{-xu}\frac{1}{2\pi i}f[n, k, u+1]  \frac{Log^{k+1}[u]}{k+1}du,
$$ 
$$
(29) = -\pi e^{-x}\sum_{k=0}^{n-1}{n\choose k}\int_0^{\infty} e^{-xu}(\frac{1}{2\pi i}\frac{d}{du}f[n, k, u+1] ) \frac{Log^{k+1}[u]}{k+1}du.
$$
To compare (27) and (24), write (24) as,
\begin{eqnarray*}
\frac{e^{-x}}{2} \int_{\frac{-3\pi}{2}}^{\frac{\pi}{2}} \Big(Log [2] +Log[\epsilon]  + i(\frac{\pi}{2} +\theta)\Big)^n  d\theta &=& \frac{e^{-x}}{2} \int_{-\pi}^{\pi} \Big(Log [2] +Log[\epsilon]  + i\varphi)^n  d\varphi \\
&=& \frac{e^{-x}}{2} \int_{-\pi}^{\pi} \sum_{k=1}^{n}{n\choose k}(Log[2]+ i\varphi)^{n-k}Log^k[\epsilon] d\varphi  \\
&+&  \frac{e^{-x}}{2} \int_{-\pi}^{\pi} (Log[2]+ i\varphi)^{n} d\varphi \qquad \qquad\qquad \quad (30)
\end{eqnarray*}
so the coefficient of $Log^k[\epsilon]$ in (24) is
$$
\frac{e^{-x}}{2}\int_{-\pi}^{\pi}{n\choose k}(Log[2]+ i\varphi)^{n-k}d\varphi  = \frac{e^{-x}}{2i}{n\choose k}\frac{f[n, k-1, 1]}{n - (k-1)}, \qquad  1\le k\le n,
$$
while the coefficient of $Log^k[\epsilon]$ in (27) is
$$
-\frac{e^{-x}}{2i}{n\choose k-1}\frac{f[n, k-1, 1]}{k}, \qquad  1\le k\le n,
$$
in which case the sum of (24) plus (27) is $0$, because,
$$
\frac{{n\choose k}}{n+1-k} - \frac{{n\choose k-1}}{k} = 0.
$$
The real part of the sum of the terms left over is the statement of Theorem 2. These terms are (30), (18), (28), (29). Listing the terms in that order, multiplying by $(-1)^n\frac{1}{2}$, gives the theorem. 
\end{proof}
\section {$\frac{d^2}{ds^2} T[s, x]\Big|_{s=.5}$}
Recall that
\begin{eqnarray*}
p_2(u) &=& 2u \\
Re\Big(\frac{1}{2\pi i}f[2, 0, u+1]\Big) &=& 2Log[u+ 2] \\
Re\Big(\frac{1}{2\pi i}f[2, 1, u+1]\Big) &=& 1 \\
Re\Big(\frac{1}{2\pi i}\frac{d}{du}f[2, 0, u+1]\Big) &=& \frac{2}{u + 2} \\
Re\Big(\frac{1}{2\pi i}\frac{d}{du}f[2, 1, u+1]\Big) &=& 0. 
\end{eqnarray*}
Therefore, from Theorem 2,
\begin{eqnarray*}
\frac{d^2}{ds^2} T[s, x]\Big|_{s = .5} &=& T[.5, x]\Big(Log^2[2] - 2 \zeta(2) \\
&-2& \int_0^{\infty}e^{-x u}\frac{Log[u + 2]}{u + 2}du -2 \int_0^{\infty}e^{-xu} \frac{Log[u]}{u + 2}du \\
&+& 2x\int_0^{\infty}e^{-x u}Log[u + 2]Log[u]du \\
&+& x\int_0^{\infty}e^{-x u}Log^2[u]du \\
&-&2 \int_0^{\infty}e^{-x u}\frac{Log[u]}{u + 2}du\Big).
\end{eqnarray*}
Two integrals combine, so that,
\begin{eqnarray*}
\frac{\ddot T}{T} &=& Log^2[2] - 2 \zeta[2] \\
&-2& \int_0^{\infty}e^{-x u}\frac{Log[u + 2]}{u + 2}du - 4 \int_0^{\infty}e^{-xu} \frac{Log[u]}{u + 2}du \\
&+& 2x\int_0^{\infty}e^{-x u}Log[u + 2]Log[u]du \\
&+& x\int_0^{\infty}e^{-x u}Log^2[u]du.
\end{eqnarray*}
With
$$
S(x) = \Big(\frac{2}{x}\Big)^s,
$$
and recalling (for $s = .5$),
$$
\frac{\dot T}{T} = \Big(\int_0^{\infty}e^{-u}\frac{du}{u+2x} +Log[x] - \dot\Gamma[1] - Log[2]\Big),
$$
some cancellation  in,
$$
\frac{\ddot K}{K} = \frac{\ddot S}{S} + \frac{\ddot \Gamma}{\Gamma} + \frac{\ddot T}{T} 
+ 2\Big(\frac{\dot S}{S} \frac{\dot \Gamma}{\Gamma}+ \frac{\dot S}{S} \frac{\dot T}{T} + \frac{\dot \Gamma}{\Gamma}\frac{\dot T}{T}\Big), 
$$
is evident. For example, with
$$
x\int_0^{\infty}e^{-x u}Log^2[u]du = \ddot\Gamma[1] - 2Log[x] \dot\Gamma[1] + Log^2[x],
$$
then among the terms of,
$$
\frac{\ddot\Gamma}{\Gamma} + 2 \frac{\dot \Gamma}{\Gamma}\frac{\dot T}{T} + \frac{\ddot T}{T}
$$
is   
$$
2\ddot\Gamma[1] - 2(\dot\Gamma[1])^2 -2\zeta[2] = 0.
$$
Other cancellation occurs, but the situation does not encourage a further look at $n\ge 3$.
\section {References}
[1] BesselK Series for the Riemann Zeta function, Timothy Redmond, Charles Ryavec, arXiv:1710.09987v2[Math.NT]. 

[2] NIST Digital Library of Mathematical Functions.
\end{document}